\documentclass[12pt]{article}
\usepackage{amsmath}
\usepackage{amssymb}
\usepackage{amscd}
\usepackage{amsthm}

\theoremstyle{plain}
\newtheorem{Thm}{Theorem}[section]

\newtheorem{Cor}[Thm]{Corollary}
\newtheorem{Lem}[Thm]{Lemma}

\theoremstyle{definition}

\newtheorem{Expl}[Thm]{Example}

\newtheorem{Rem}[Thm]{Remark}

\numberwithin{equation}{section}

\title{Euivalences of derived catgories of sheaves on smooth stacks}
\author{Yujiro Kawamata}

\begin{document}

\maketitle

\begin{abstract}
We extend Orlov's representability theorem on the equivalence of 
derived categories of sheaves to the case of smooth stacks associated to 
normal projective varieties with only quotient singularities.
\end{abstract}


\section{Introduction}

Orlov's representability theorem \cite{Orlov} is of fundamental importance
in the study of derived categories of sheaves on smooth projective varieties.
It says that any equivalence of derived categories for such 
varieties is representable as a Fourier-Mukai functor.
The purpose of this paper is to extend this theorem to the case of projective 
orbifolds.

To any normal projective variety $X$ which has only quotient singularities,
we can attach naturally a smooth stack (or an orbifold) $\mathcal{X}$.
Many results on smooth varieties are expected to extend to such 
stacks. 
For example, any coherent sheaf on a smooth stack has a finite locally
free resolution provided that the stack has a trivial stabilizer group at the
generic point and the coarse moduli space is a separated scheme 
(\cite{Totaro}).
On the other hand, if we allow quotient singularities on varieties, 
then there are a lot
more possibility of interesting examples such as flops than the case of 
smooth varieties.
Note that there are more sheaves on the stack $\mathcal{X}$ 
than the underlying variety $X$, and as was shown in \cite{Francia},
the correct category which should be considered is that of the stacky 
sheaves, not the ordinary sheaves.

For a smooth stack $\mathcal{X}$, we denote by $D^b(\text{Coh}(\mathcal{X}))$,
or sometimes simply $D^b(\mathcal{X})$,
the derived category of bounded complexes of coherent sheaves.
If $e \in D^b(\text{Coh}(\mathcal{X} \times \mathcal{Y}))$ 
is an object on the product stack, the {\it integral functor} 
\[
\Phi^e_{\mathcal{X} \to \mathcal{Y}}: D^b(\text{Coh}(\mathcal{X}))
\to D^b(\text{Coh}(\mathcal{Y}))
\]
between such categories is defined by
\[
\Phi^e_{\mathcal{X} \to \mathcal{Y}}(a) = p_{2*}(e \otimes p_1^*a)
\]
for $a \in D^b(\text{Coh}(\mathcal{X}))$, where
$p_1^*$ and $\otimes$ are the left derived functors and
$p_{2*}$ is the right derived functor.
An integral functor is called a {\it Fourier-Mukai functor} if it is an
equivalence.

The main result of this paper is 
the following extension of \cite{Orlov}~Theorem~2.2:

\begin{Thm}\label{main}
Let $X$ and $Y$ be normal projective varieties with only quotient singularities
and let $\mathcal{X}$ and $\mathcal{Y}$ be smooth stacks naturally associated 
to them.
Let 
\[
F: D^b(\text{Coh}(\mathcal{X})) \to D^b(\text{Coh}(\mathcal{Y}))
\]
be an exact functor which is fully faithful and has a left adjoint functor.
Then there exist an object 
$e \in D^b(\text{Coh}(\mathcal{X} \times \mathcal{Y}))$
and an isomorphism of functors
\[
F \cong \Phi^e_{\mathcal{X} \to \mathcal{Y}}
\]
Moreover, the object $e$ is uniquely determined up to isomorphism.
\end{Thm}

In \S 2 we review the way to attach a single object to a complex of 
objects.  
There are two different ways of construction, the left and right convolutions,
which are both used in the proof of Theorem~\ref{main}.

In \S 3 we construct an infinite left resolution of the structure sheaf 
of the diagonal of the self product for any projective variety
(Theorem~\ref{diagonal1}).
The construction is based on a theorem of Backelin \cite{Backelin} 
which says that the
Veronese subring of the homogeneous coordinate ring 
of sufficiently high degree becomes a Koszul algebra.

In \S 4 we extend the construction of \S 3 to smooth stacks by taking the
invariant part under a group action (Theorem~\ref{diagonal2}).
We also prove that for any coherent sheaf on our stack, there exists 
a surjective morphism from a locally free sheaf of special type
(Theorem~\ref{locallyfree}).
This gives an alternative proof of Totaro's theorem (\cite{Totaro}) 
for our stack.
As an example, we describe the resolution of the diagonal explicitly in the
case of weighted projective spaces using Canonaco's formula \cite{Canonaco}
in \S 5.

We prove the main theorem in \S 6.
The proof is partly parallel to that of \cite{Orlov}, but we also need
different arguments because there are more sheaves on the stacks.
Indeed, we need both left and right convolutions, an infinite resolution
of the diagonal and a special locally free sheaf, which are explained 
in previous sections.

\S 7 is devoted to similar applications to those appearing in \cite{DK}. 
We consider the consequences of \lq \lq $D$-equivalence" in 
Theorem~\ref{application1}, and for the group of autoequivalences in 
Theorem~\ref{application2}.

We work over $\mathbb{C}$, but most of our results hold for arbitrary fields. 

The author would like to thank Keiji Oguiso for pointing out the reference
\cite{FH}.
He would also like to thank the referee for many valuable comments.


\section{Convolution}

We sometimes have to deal with a complex of objects in a derived category.
But we cannot associate a single object to it in a similar way as a single 
complex associated to a double complex, because  
quasi-isomorphisms are inverted in the derived category.

We recall the theory of convolutions from \cite{Orlov}.
Let $\mathcal{D}$ be a triangulated category, and let 
\begin{equation}\label{complex1}
\begin{CD}
0 @>>> a_m @>{d_m}>> a_{m-1} @>{d_{m-1}}>> \dots @>>> a_1 @>{d_1}>> a_0 @>>> 0
\end{CD}
\end{equation}
be a bounded complex of objects in $\mathcal{D}$, that is, we assume that
$d_id_{i+1} = 0$ for $1 \le i < m$.
Assume the condition 
\begin{equation}\label{complex1cond}
\text{Hom}_{\mathcal{D}}(a_p[r], a_q) = 0 \text{ for any }
p > q, r > 0.
\end{equation}
Then the {\it right convolution} is defined to be a morphism
$d_0: a_0 \to a$ to an object in $\mathcal{D}$ which is uniquely constructed 
inductively on the length $m$ in the following way.
(We would like to call this \lq \lq right" convolution 
instead of \lq \lq left" as in \cite{Orlov} 
because the direction of the convolution is to the right.)

If $m = 0$, then $a = a_0$ and $d_0 = \text{Id}$.  
If $m \ge 1$, 
then let $a'_{m-1}$ be the cone of the morphism $a_m \to a_{m-1}$:
\[
\begin{CD}
a_m @>>> a_{m-1} @>{j_{m-1}}>> a'_{m-1} @>>> a_m[1]
\end{CD}
\]
is a distinguished triangle.
We have an exact sequence
\[
\begin{split}
&\text{Hom}_{\mathcal{D}}(a_m[1], a_{m-2}) \to 
\text{Hom}_{\mathcal{D}}(a'_{m-1}, a_{m-2}) \\ &\to 
\text{Hom}_{\mathcal{D}}(a_{m-1}, a_{m-2}) \to 
\text{Hom}_{\mathcal{D}}(a_m, a_{m-2}).
\end{split}
\]
Since the first term as well as the image of $d_{m-1}$ 
in the last term vanishes, there exists a unique morphism
$d'_{m-1}: a'_{m-1} \to a_{m-2}$ such that 
$d'_{m-1}j_{m-1} = d_{m-1}$.
Moreover, we deduce that $d_{m-2}d'_{m-1} = 0$ by a similar diagram chasing.

Thus we obtain a new complex
\begin{equation}\label{complex2}
\begin{CD}
0 @>>> a'_{m-1} @>{d'_{m-1}}>> a_{m-2} @>{d_{m-2}}>> 
\dots @>>> a_1 @>{d_1}>> a_0 @>>> 0
\end{CD}
\end{equation}
for which we can check condition $(\ref{complex1cond})$.

The right convolution $d_0: a_0 \to a$ is obtained as that 
of the complex $(\ref{complex2})$.
Note that the object $a_0$ is unchanged during the above process.
Thus we obtained:

\begin{Lem}[\cite{Orlov}~Lemma~1.5]\label{complex1lem}
Let $(\ref{complex1})$ be a complex of objects in a triangulated category
$\mathcal{D}$ which satisfies condition $(\ref{complex1cond})$.
Then there exists a right convolution
$d_0: a_0 \to a$ which is uniquely determined up to isomorphism.
\end{Lem}

\begin{Expl}
(1) If all the $a_p$ are sheaves, then the right convolution is 
nothing but the complex itself with an obvious morphism $a_0 \to a_{\bullet}$.
Indeed, we can take the intermediate objects $a'_p$ to be the complex
\[
a'_p = \{a_m \to a_{m-1} \to \dots \to a_p\}
\]
where $a_p$ is at degree $0$.

(2) If $d_0: a_0 \to a$ is the right convolution of the 
complex $(\ref{complex1})$ and $F$ is a fully faithful exact functor, 
then $F(d_0): F(a_0) \to F(a)$ 
is the right convolution of the complex
\[
0 \to F(a_m) \to F(a_{m-1}) \to \dots \to F(a_1) \to F(a_0) \to 0.
\]
Indeed, the assumption is preserved by $F$ and 
the objects of intermediate steps are given by $F(a'_p)$.
\end{Expl}

\begin{Lem}[cf. \cite{Orlov}~Lemma~1.6]\label{complexhom1lem}
Let 
\begin{equation}\label{complexhom1}
\begin{CD}
0 @>>> a_m @>{d_m}>> a_{m-1} @>{d_{m-1}}>> \dots @>>> a_1 @>{d_1}>> a_0 @>>> 0
\\
@.     @V{f_m}VV     @V{f_{m-1}}VV         @.         @V{f_1}VV     @V{f_0}VV
@. \\
0 @>>> b_m @>{e_m}>> b_{m-1} @>{e_{m-1}}>> \dots @>>> b_1 @>{e_1}>> b_0 @>>> 0
\end{CD}
\end{equation}
be a morphism between complexes which satisfy condition
$(\ref{complex1cond})$.
Let $d_0: a_0 \to a$ and $e_0: b_0 \to b$ 
be the right convolutions of these complexes, and $\epsilon: b \to b'$ 
a morphism.
Assume in addition that 
\begin{equation}\label{complexhom1cond1}
\text{Hom}_{\mathcal{D}}(a_p[r], b_q) = 0 \text{ for any } p > q, r > 0.
\end{equation}
Then there exists a morphism
$f: a \to b'$ in $\mathcal{D}$ 
satisfying the commutativity $fd_0 = \epsilon e_0f_0$:
\begin{equation}\label{morphcom}
\begin{CD}
a_0 @>{d_0}>> a @>=>> a\\
@V{f_0}VV @. @VfVV \\
b_0 @>{e_0}>> b @>{\epsilon}>> b'.
\end{CD}
\end{equation}
Moreover, this morphism is unique with this property if
\begin{equation}\label{complexhom1cond2}
\text{Hom}_{\mathcal{D}}(a_p[r], b') = 0 \text{ for any } p, r > 0.
\end{equation}
\end{Lem}

\begin{proof}
In order to find $f$, we proceed by induction on the length $m$.
If $m = 0$, then $f = \epsilon f_0$.  
If $m \ge 1$, then we construct $a'_{m-1}$ and $b'_{m-1}$ as before.
The morphisms $f_m$ and $f_{m-1}$ induce a morphism $f'_{m-1}: a'_{m-1} \to 
b'_{m-1}$ between the cones.
We note that the existence of the morphism $f'_{m-1}$ is guaranteed by the 
axiom of triangulated categories such that the following diagram
\[
\begin{CD}
a_m @>>> a_{m-1} @>{j_{m-1}}>> a'_{m-1} @>>> a_m[1] \\
@V{f_m}VV @V{f_{m-1}}VV        @V{f'_{m-1}}VV @V{f_m[1]}VV \\
b_m @>>> b_{m-1} @>{k_{m-1}}>> b'_{m-1} @>>> b_m[1].
\end{CD}
\]
is commutative, but it is not uniquely determined in general.

We have
\[
e'_{m-1}f'_{m-1}j_{m-1} = e'_{m-1}k_{m-1}f_{m-1}
= e_{m-1}f_{m-1} = f_{m-2}d_{m-1} = f_{m-2}d'_{m-1}j_{m-1}.
\]
Hence $e'_{m-1}f'_{m-1} = f_{m-2}d'_{m-1}$ by condition 
$(\ref{complexhom1cond1})$.
By diagram chasing, we can check that condition 
$(\ref{complexhom1cond1})$ is satisfied by the new morphism of complexes
\begin{equation}\label{complexhom2}
\begin{CD}
0 @>>> a'_{m-1} @>{d'_{m-1}}>> a_{m-2} @>{d_{m-2}}>> \dots @>>> 
a_1 @>{d_1}>> a_0 @>>> 0 \\
@.     @V{f'_{m-1}}VV     @V{f_{m-2}}VV   @.   @V{f_1}VV     @V{f_0}VV   @. \\
0 @>>> b'_{m-1} @>{e'_{m-1}}>> b_{m-2} @>{e_{m-2}}>> \dots @>>> 
b_1 @>{e_1}>> b_0 @>>> 0
\end{CD}
\end{equation}
from which we obtain the morphism $f: a \to b'$ by the induction hypothesis.

We note that $f_0: a_0 \to b_0$ remains unchanged in the above process.
By condition $(\ref{complexhom1cond2})$, 
$f$ is uniquely determined by the commutativity of $(\ref{morphcom})$.
\end{proof}


In the dual way, we define the left convolutions of complexes of objects.
We consider a bounded complex of objects $(\ref{complex1})$
which satisfies condition $(\ref{complex1cond})$.
The {\it left convolution} is a morphism $d: a \to a_m$ from 
an object in $\mathcal{D}$ which is uniquely constructed 
inductively on the length $m$ in the following way.

If $m = 0$, then $a = a_0$ and $d = \text{Id}$.  
If $m \ge 1$, 
then let $a'_1[1]$ be the cone of the morphism $a_1 \to a_0$:
\[
\begin{CD}
a'_1 @>{j_1}>> a_1 @>{d_1}>> a_0 @>>> a'_1[1]
\end{CD}
\]
is a distinguished triangle.
There exists a unique morphism
$d'_2: a_2 \to a'_1$ such that $j_1d'_2 = d_2$ and $d'_2d_3 = 0$ as before.

Thus we obtain a new complex
\begin{equation}\label{complex3}
\begin{CD}
0 @>>> a_m @>{d_m}>> a_{m-1} @>{d_{m-1}}>> 
\dots @>>> a_2 @>{d'_2}>> a'_1 @>>> 0
\end{CD}
\end{equation}
which satisfies condition $(\ref{complex1cond})$.
The left convolution $d: a \to a_m$ of the complex $(\ref{complex1})$ is 
defined to be that of the new complex $(\ref{complex3})$. 

\begin{Lem}
Let $(\ref{complex1})$ be a complex of objects in a triangulated category
$\mathcal{D}$ which satisfies condition $(\ref{complex1cond})$.
Then there exists a left convolution
$d: a \to a_m$ which is uniquely determined up to isomorphism.
\end{Lem}

\begin{Expl}
(1) If all the $a_p$ are sheaves, then the left convolution is 
nothing but the complex itself with an obvious morphism 
$a_{\bullet}[-m] \to a_m$, where the term $a_m$ of $a_{\bullet}[-m]$ 
is at degree $0$.
Indeed, we can take inductively
\[
a'_p = \{a_p \to a_{p-1} \to \dots \to a_1 \to a_0\}
\]
where $a_p$ is at degree $0$.

(2) If $d: a \to a_m$ is the left convolution of the complex 
$(\ref{complex1})$ 
and $F$ is a fully faithful exact functor, 
then $F(d): F(a) \to F(a_m)$ is the left convolution of the complex
\[
0 \to F(a_m) \to F(a_{m-1}) \to \dots \to F(a_1) \to F(a_0) \to 0.
\]
\end{Expl}

\begin{Lem}
Let $(\ref{complexhom1})$ be a morphism between complexes 
which satisfy condition $(\ref{complex1cond})$.
Let $d: a \to a_m$ and $e: b \to b_m$ be the left convolutions of these 
complexes, and $\epsilon: a' \to a$ a morphism.
Assume in addition that 
\[
\text{Hom}_{\mathcal{D}}(a_p[r], b_q) = 0 \text{ for any }
p > q, r > 0
\]
Then there exists a morphism
$f: a' \to b$ in $\mathcal{D}$ 
satisfying the commutativity $ef = f_md \epsilon$:
\[
\begin{CD}
a' @>{\epsilon}>> a @>d>> a_m \\
@VfVV @. @V{f_m}VV \\
b @>=>> b @>e>> b_m.
\end{CD}
\]
Moreover, this morphism is unique with this property if
\[
\text{Hom}_{\mathcal{D}}(a'[r], b_q) = 0 \text{ for any } q, r > 0.
\]
\end{Lem}

\begin{proof}
The proof is similar to that of Lemma~\ref{complexhom1lem}.
\end{proof}


\section{Resolution of the diagonal}

For a projective space $\mathbb{P}$, 
there is a standard resolution of the structure sheaf 
of the diagonal subvariety $\mathcal{O}_{\Delta \mathbb{P}}$
in the self-product 
$\Delta \mathbb{P} \subset \mathbb{P} \times \mathbb{P}$ by
locally free sheaves of K\"unneth type (\cite{Beilinson}, see also
$(\ref{Beilinson})$).
We extend this in a weaker sense for general projective varieties.

Let $X$ be a projective algebraic variety, $L$ an ample line bundle,
and let 
\[
A = \bigoplus_{m=0}^{\infty} A_m = \bigoplus_{m=0}^{\infty} H^0(X, mL)
\]
be the homogeneous coordinate ring.
We define vector spaces $B_m$ ($m \ge 0$) by $B_0 = A_0$, $B_1 = A_1$ and
\[
B_m = \text{Ker}(B_{m-1} \otimes_{\mathbb{C}} A_1 \to 
B_{m-2} \otimes_{\mathbb{C}} A_2)
\]
for $m \ge 2$, where the homomorphism is induced from the multiplication in  
$A$.
We set $A_m = B_m = 0$ for $m < 0$.

The graded $\mathbb{C}$-algebra $A$ is said to be a {\it Koszul algebra} 
if the sequence of natural homomorphisms 
\begin{equation}\label{seq-module}
\begin{split}
&\dots \to B_m \otimes_{\mathbb{C}} A(-m) \to 
B_{m-1} \otimes_{\mathbb{C}} A(-m+1) \to 
\dots \\
&\to B_1 \otimes_{\mathbb{C}} A(-1) \to A \to \mathbb{C} \to 0
\end{split}
\end{equation}
is exact, where the shifted module $A(j)$ is defined by 
$A(j)_k = A_{j+k}$.  
In other words, $A$ is Koszul if the minimal $A$-free
resolution of its residue field consists of homomorphisms of degree $1$.

By \cite{Backelin}~Theorem~2, the subring 
$A^{(d)} = \bigoplus_{m=0}^{\infty} A_{dm}$ of $A$ is a Koszul algebra
for a sufficiently large integer $d$, i.e.,
$A$ becomes Koszul if we replace $L$ by its tensor power $L^d$.
We assume that $L$ is already replaced so that $A$ is Koszul in the following.

We define sheaves $R_m$ ($m \ge 0$) on $X$ by $R_0 = \mathcal{O}_X$ and 
\[
R_m = \text{Ker}(B_m \otimes_{\mathbb{C}} \mathcal{O}_X \to 
B_{m-1} \otimes_{\mathbb{C}} L)
\]
for $m \ge 1$, where the homomorphism is induced from the natural 
homomorphisms $B_m \to B_{m-1} \otimes_{\mathbb{C}} A_1$ and 
$A_1 \otimes_{\mathbb{C}} \mathcal{O}_X \to L$.
We set $R_m = 0$ for $m < 0$.

If we take the sequence of associated sheaves to
the tensor product of $(\ref{seq-module})$ with $A(m)$,
we obtain an exact sequence
\begin{equation}\label{seq-sheaf}
0 \to R_m \to B_m \otimes_{\mathbb{C}} \mathcal{O}_X
\to B_{m-1} \otimes_{\mathbb{C}} L \to 
\dots \to B_1 \otimes_{\mathbb{C}} L^{m-1} \to L^m \to 0.
\end{equation}

\begin{Lem}\label{AR}
There is an exact sequence
\[
\begin{split}
&0 \to A_0 \otimes_{\mathbb{C}} R_m \to A_1 \otimes_{\mathbb{C}} R_{m-1} 
\to \\
&\dots \to 
A_{m-1} \otimes_{\mathbb{C}} R_1 \to A_m \otimes_{\mathbb{C}} 
R_0 \to L^m \to 0.
\end{split}
\]
\end{Lem}

\begin{proof}
We consider a double complex of sheaves
\[
C^{i,j} = \begin{cases}
A_i \otimes B_{m-i-j} \otimes L^j &\text{ for } i, j, m - i - j \ge 0 \\
0 &\text{ otherwise }
\end{cases}
\]
where the differentials
\[
\begin{split}
&d_1^{i,j}: A_i \otimes B_{m-i-j} \otimes L^j
\to A_{i+1} \otimes B_{m-i-j-1} \otimes L^j \\
&d_2^{i,j}: A_i \otimes B_{m-i-j} \otimes L^j
\to A_i \otimes B_{m-i-j-1} \otimes L^{j+1}
\end{split}
\]
are induced from the homomorphisms in $(\ref{seq-module})$ 
and $(\ref{seq-sheaf})$.
The condition $d_1^{i,j+1}d_2^{i,j} = \pm d_2^{i+1,j}d_1^{i,j}$
is satisfied, because both are obtained as a composition of natural 
homomorphisms
\[
\begin{split}
&A_i \otimes B_{m-i-j} \otimes L^j \to 
A_i \otimes B_{m-i-j-2} \otimes A_1 \otimes A_1 \otimes L^j \\ &\to  
A_{i+1} \otimes B_{m-i-j-2} \otimes L^{j+1}.
\end{split}
\]

If we take the cohomologies of $C^{i,j}$ with respect to the differential 
$d_1$ first, 
then we obtain $0$'s except $L^m$ at $(i,j)=(0,m)$.
On the other hand, if we take the cohomologies with respect to $d_2$ first, 
then we obtain $A_i \otimes R_{m-i}$ at $j = 0$ and $0$ elsewhere
by $(\ref{seq-sheaf})$.
Therefore, we obtain our exact sequence.
\end{proof}

Now we consider sheaves on the product $X \times X$.
By Lemma~\ref{AR}, we obtain a homomorphism 
$L^{-m} \boxtimes R_m \to L^{-m+1} \boxtimes R_{m-1}$
as a composition of the following homomorphisms
\begin{equation}\label{LRhom}
L^{-m} \boxtimes R_m \to A_1 \otimes_{\mathbb{C}} (L^{-m} \boxtimes R_{m-1})
\to L^{-m+1} \boxtimes R_{m-1}.
\end{equation}

\begin{Thm}\label{diagonal1}
Let $X$ be a projective algebraic variety, $L$ an ample line bundle, and
$\Delta X \subset X \times X$ the diagonal subvariety of the direct product.
If $L$ is replaced by a sufficiently high power $L^d$, then 
the complex of sheaves on $X \times X$ 
\[
\begin{split}
&\dots \to L^{-m} \boxtimes R_m \to L^{-m+1} \boxtimes R_{m-1} \to 
\dots \\ &\to L^{-1} \boxtimes R_1 \to \mathcal{O}_X \boxtimes \mathcal{O}_X
\to \mathcal{O}_{\Delta X} \to 0
\end{split}
\]
is exact.
\end{Thm}

\begin{proof}
Let $D^{\bullet}$ denote the given complex, where we set 
$D^m = L^{m+1} \boxtimes R_{-m-1}$ for $m \le -1$ and 
$D^0 = \mathcal{O}_{\Delta X}$.
Assume that $H^{m_0}(D^{\bullet}) \ne 0$ for some non-positive integer $m_0$.
We take a sufficiently large integer $k$ such that 
\[
\begin{split}
&R^pp_{2*}(H^q(D^{\bullet}) \otimes p_1^*L^k) = 0 \text{ for } 
p > 0, q \ge m_0 - \dim X \\
&H^p(X, L^{k+q}) = 0 \text{ for } p > 0, q \ge m_0 - \dim X \\
&R^0p_{2*}(H^{m_0}(D^{\bullet}) \otimes p_1^*L^k) \ne 0.
\end{split}
\]
We consider a spectral sequence
\[
^{'}E_2^{p,q} = R^pp_{2*}(H^q(D^{\bullet}) \otimes p_1^*L^k) \Rightarrow 
R^{p+q}p_{2*}(D^{\bullet} \otimes p_1^*L^k).
\]
By the assumption, we have 
$^{'}E_2^{p,q} = 0$ for $p > 0$ and $q \ge m_0 - \dim X$
while $^{'}E_2^{0,m_0} \ne 0$, 
hence $R^{m_0}p_{2*}(D^{\bullet} \otimes p_1^*L^k) \ne 0$.

On the other hand, we consider another spectral sequence
\[
{}^{''}E_1^{p,q} = R^qp_{2*}(D^p \otimes p_1^*L^k) \Rightarrow 
R^{p+q}p_{2*}(D^{\bullet} \otimes p_1^*L^k).
\]
We have
\[
{}^{''}E_1^{p,q} = 
\begin{cases} L^k &\text{ for } p = 0 \text{ and } q = 0 \\
H^q(X, L^{k+p+1}) \otimes R_{-p-1} &\text{ for } p \le -1 \\
0 &\text{ otherwise.} 
\end{cases}
\]
Thus ${}^{''}E_1^{p,0} = A_{k+p+1} \otimes R_{-p-1}$ for $p \le -1$ and 
${}^{''}E_1^{p,q} = 0$ for $p + 1 \ge m_0 - \dim X$ and $q > 0$.
By Lemma~\ref{AR}, we obtain that ${}^{''}E_2^{p,q} = 0$ for 
$p + q = m_0$, a contradiction.
\end{proof}

\begin{Cor}\label{Id1}
Assume the conditions of Theorem~\ref{diagonal1}.
Define an object $e \in D^-(X \times X)$ by
\[
\begin{split}
e = \{&\dots \to L^{-m} \boxtimes R_m \to L^{-m+1} \boxtimes R_{m-1} \to 
\dots \\ &\to L^{-1} \boxtimes R_1 \to \mathcal{O}_X \boxtimes \mathcal{O}_X\}.
\end{split}
\]
Then for any object $a \in D^-(X)$,
there are isomorphisms in $D^-(X)$
\[
a \cong p_{1*}(e \otimes p_2^*a) \cong p_{2*}(e \otimes p_1^*a)
\]
where $p_{1*}$ and $p_{2*}$ are right derived functors and 
$\otimes$, $p_1^*$ and $p_2^*$ are left derived functors. 
\end{Cor}

\begin{proof}
We have $p_{1*}(\mathcal{O}_{\Delta X} \otimes p_2^*a) \cong 
p_{2*}(\mathcal{O}_{\Delta X} \otimes p_1^*a) \cong a$.
\end{proof}


\section{Smooth stack}

We refer to \cite{Mumford} and \cite{Gomez} for readable account on 
Deligne-Mumford stacks.

Let $X$ be a normal projective variety with only quotient singularities.
Then there exists an etale covering $\{U_i\}$ of $X$ with finite Galois 
coverings $\sigma_i: \mathcal{U}_i \to U_i$ from smooth varieties 
$\mathcal{U}_i$ which are etale in codimension $1$.
The data $\{U_i, \mathcal{U}_i, \sigma_i\}$ defines a 
{\it smooth stack} $\mathcal{X}$.
Let $\sigma: \mathcal{X} \to X$ denote the natural morphism.

Let $X^{\#}$ be the normalization of $X$ in a common Galois extension 
of the fields $\mathbb{C}(\mathcal{U}_i)$.
Then we have a morphism $\pi: X^{\#} \to \mathcal{X}$.
Denote $G = \text{Gal}(X^{\#}/X)$ and $G_i = \text{Gal}(\mathcal{U}_i/U_i)$.
We have the following diagram
\[
\begin{CD}
U^{\#}_i @>{\pi_i}>> \mathcal{U}_i @>{\sigma_i}>> U_i \\
@V{\text{open}}VV @V{\text{etale}}VV @V{\text{etale}}VV \\
X^{\#} @>{\pi}>> \mathcal{X} @>{\sigma}>> X.
\end{CD}
\]
where $U^{\#}_i$ is an open subset of $X^{\#}$ such that the induced morphism
$\pi_i: U^{\#}_i \to \mathcal{U}_i$ is finite.

Let $a$ be a sheaf on $X^{\#}$. 
Then the {\it direct image sheaf} $\pi_*a$ on $\mathcal{X}$ is given by the 
following collection of $G_i$-equivariant sheaves $(\pi_*a)_i$ 
on the $\mathcal{U}_i$.
Let $X^{\#}_i = (X^{\#} \times_X \mathcal{U}_i)^{\nu}$ 
be the normalization of the fiber product
with natural morphisms $p_1: X^{\#}_i \to X^{\#}$
and $p_2: X^{\#}_i \to \mathcal{U}_i$.
Then 
\[
(\pi_*a)_i = p_{2*}p_1^*a.
\]
We note that $p_1$ is etale.
Indeed, the normalization of the fiber product
$\mathcal{U}_{ij} = (\mathcal{U}_i \times_X \mathcal{U}_j)^{\nu}$ 
is etale over $\mathcal{U}_j$, and 
$X^{\#}_i \cap p_1^{-1}(U_j^{\#}) = U_j^{\#} \times_{\mathcal{U}_j}
\mathcal{U}_{ij}$.

If $b$ is a sheaf on $\mathcal{X}$ given by the collection of 
$G_i$-equivariant sheaves $b_i$ on the $\mathcal{U}_i$, 
then the {\it inverse image sheaf} 
$\pi^*b$ on $X^{\#}$ is obtained as the etale descent of the 
collection of sheaves $p_2^*b_i$ on the $X^{\#}_i$.

Since 
\[
(\pi_*\pi^*b)_i = p_{2*}p_2^*b_i
\]
the Galois group $G$ acts on the sheaf $\pi_*\pi^*b$, and we have
\[
(\pi_*\pi^*b)^G \cong b.
\] 

Let $G = \text{Gal}(X^{\#}/X) \cong \Delta G \subset G \times G$ 
act diagonally on the product $X^{\#} \times X^{\#}$. 

\begin{Thm}\label{diagonal2}
Let $X$ be a normal projective variety with only quotient singularities, 
$\sigma: \mathcal{X} \to X$ the natural morphism from the 
associated smooth stack, $L$ an ample line bundle, $\mathcal{L} = 
\sigma^*L$, 
and $\pi: X^{\#} \to \mathcal{X}$ 
a finite Galois morphism from a projective variety as above.
Define the sheaves $R_m$ on $X^{\#}$ from $\pi^*\mathcal{L}$ as before.
Define an object $e^- \in D^-(\text{Coh}(\mathcal{X} \times \mathcal{X}))$ by
\begin{equation}\label{diagonaleq2}
\begin{split}
e^- = \{&\dots \to (\pi_*\pi^*\mathcal{L}^{-m} \boxtimes \pi_*R_m)^{\Delta G} 
\to (\pi_*\pi^*\mathcal{L}^{-m+1} \boxtimes \pi_*R_{m-1})^{\Delta G} \to 
\dots \\ &\to (\pi_*\pi^*\mathcal{L}^{-1} \boxtimes \pi_*R_1)^{\Delta G} \to 
(\pi_*\pi^*\mathcal{O}_{\mathcal{X}} 
\boxtimes \pi_*\pi^*\mathcal{O}_{\mathcal{X}})^{\Delta G}
\}.
\end{split}
\end{equation}
If $L$ is replaced by a sufficiently high power $L^d$, 
then there are isomorphisms in $D^-(\text{Coh}(\mathcal{X}))$ 
\[
a \cong p_{1*}(e^- \otimes p_2^*a) \cong p_{2*}(e^- \otimes p_1^*a).
\]
for any object $a \in D^-(\text{Coh}(\mathcal{X}))$
\end{Thm}

\begin{proof}
Let 
\[
\begin{split}
e^{\#} = \{&\dots \to \pi^*\mathcal{L}^{-m} \boxtimes R_m \to 
\pi^*\mathcal{L}^{-m+1} \boxtimes R_{m-1} 
\to 
\dots \\ &\to \pi^*\mathcal{L}^{-1} \boxtimes R_1 \to 
\pi^*\mathcal{O}_{\mathcal{X}} \boxtimes \pi^*\mathcal{O}_{\mathcal{X}}\}
\in D^-(X^{\#} \times X^{\#}).
\end{split}
\]
The homomorphisms are equivariant with respect 
to the action of $G$.
Indeed, we have a commutative diagram of injective homomorphisms
\[
\begin{CD}
R_m @>>> B_m \otimes \mathcal{O}_{X^{\#}} \\
@VVV @VVV \\
A_1 \otimes R_{m-1} @>>> A_1 \otimes B_{m-1} \otimes \mathcal{O}_{X^{\#}}.
\end{CD}
\]   
Since the multiplication $A \otimes_{\mathbb{C}} A \to A$ is $G$-equivariant
with respect 
to the diagonal action, so is the right vertical arrow, hence the left.
The homomorphism $A_1 \otimes \pi^*\mathcal{L}^{-m} \to 
\pi^*\mathcal{L}^{-m+1}$ is similarly $G$-equivariant.
Therefore, the homomorphisms in $(\ref{LRhom})$ are equivariant 
with respect to the diagonal $G$-actions.

We have $G$-equivariant isomorphisms
\[
\pi^*a \cong p_{1*}(e^{\#} \otimes p_2^*\pi^*a)
\cong p_{2*}(e^{\#} \otimes p_1^*\pi^*a)
\]
by Corollary~\ref{Id1}.
Thus 
\[
\begin{split}
&a \cong (\pi_*\pi^*a)^G \cong (\pi_*p_{1*}(e^{\#} \otimes p_2^*\pi^*a))^G
\cong (p_{1*}(\pi \times \pi)_*(e^{\#} \otimes p_2^*\pi^*a))^G \\
&\cong (p_{1*}((\pi \times \pi)_*e^{\#} \otimes p_2^*a))^G
\cong p_{1*}(((\pi \times \pi)_*e^{\#})^G \otimes p_2^*a)
\cong p_{1*}(e^- \otimes p_2^*a).
\end{split}
\]
The second isomorphism is obtained similarly.
\end{proof}

Since there are more sheaves on the stack $\mathcal{X}$ than the 
underlying variety $X$, 
the existence of a surjective homomorphism from 
a locally free sheaf as in the following is non-trivial.
This gives an alternative proof of a theorem by Totaro \cite{Totaro} 
which says that there exists a finite locally free resolution for any
coherent sheaf on $\mathcal{X}$.

\begin{Thm}\label{locallyfree}
Let $X$ be a normal quasi-projective variety with only quotient singularities,
$L$ a very ample invertible sheaf, 
$\sigma: \mathcal{X} \to X$ the morphism from the associated smooth stack, 
and $\mathcal{L} = \sigma^*L$.
Then there exists a locally free sheaf $\mathcal{A}_0$ on $\mathcal{X}$ 
such that, 
for any coherent sheaf $\mathcal{C}$ on $\mathcal{X}$, 
there exist positive integers 
$k,l$ and a surjective homomorphism 
$(\mathcal{A}_0 \otimes \mathcal{L}^{-k})^{\oplus l} \to \mathcal{C}$.
\end{Thm}

\begin{proof}
Let $T$ be the tangent sheaf of $\mathcal{X}$.
There exists a positive integer $j_0$ such that 
any irreducible representation of the stabilizer group $G_x$ of any 
point $x \in \mathcal{X}$ appears in the representation
$(\bigoplus_{j=0}^{j_0} T^{\otimes j}) \otimes \mathcal{O}_x$ of $G_x$
(\cite{FH}~Problem~2.37).

We set $\mathcal{A}_0 = (\bigoplus_{j=0}^{j_0} T^{\otimes j})$.
Let $\mathcal{C}$ be any non-zero coherent sheaf on $\mathcal{X}$, and 
$x \in \mathcal{X}$ a point in its support. 
Then there is a non-trivial homomorphism 
$\mathcal{A}_0 \otimes \mathcal{O}_x \to \mathcal{C} \otimes \mathcal{O}_x$.
Hence the sheaf
$\sigma_*\mathcal{H}om(\mathcal{A}_0, \mathcal{C})$ on $X$ is non-trivial, and
there exists a positive integer $k$ such that 
$\text{Hom}(\mathcal{A}_0 \otimes \mathcal{L}^{-k}, \mathcal{C}) \ne 0$.  
Let 
\[
\mathcal{C}' = \text{Im}(\text{Hom}(\mathcal{A}_0 \otimes \mathcal{L}^{-k}, 
\mathcal{C}) 
\otimes \mathcal{A}_0 \otimes \mathcal{L}^{-k}, \mathcal{C}).
\]
If $\mathcal{C}/\mathcal{C}' \ne 0$, 
then there exists an integer $k' \ge k$ such that 
$\text{Hom}(\mathcal{A}_0 \otimes \mathcal{L}^{-k'}, \mathcal{C}/\mathcal{C}') 
\ne 0$ and 
$\text{Hom}(\mathcal{A}_0 \otimes \mathcal{L}^{-k'}, \mathcal{C}'[1]) = 0$.
By Neotherian induction, we obtain our result.
\end{proof}


\section{An example: weighted projective space}

In this section, we consider the special case where 
$X$ is a weighted projective space as an example of the general construction.
The results will not be used in later sections.

Let $(a_0, \dots, a_n)$ be a sequence of positive integers which is 
{\it well formed}, 
that is, any subset consisting of $n$ integers from these $n+1$ is coprime.
Let a finite group $G = \mu_{a_0} \times \dots \times \mu_{a_n}$ act on 
$\mathbb{P}^{\#} = \mathbb{P}^n$ with homogeneous coordinates 
$(x^{\#}_0, \dots, x^{\#}_n)$ diagonally
so that $\mathbb{P} = \mathbb{P}^{\#}/G \cong 
\mathbb{P}(a_0, \dots, a_n)$.
Let $\pi': \mathbb{P}^{\#} \to \mathbb{P}$ be the projection.

The weighted projective space $\mathbb{P}$ has only quotient singularities.  
Let $\mathcal{P}$ be the smooth stack associated to $\mathbb{P}$, 
and $\pi: \mathbb{P}^{\#} \to \mathcal{P}$ the 
morphism induced by $\pi'$.
The stack $\mathcal{P}$ is described as follows.
Let $D^{\#}_i = \text{div}(x^{\#}_i)$,
$U^{\#}_i = \mathbb{P}^{\#} \setminus D^{\#}_i$, and 
$U_i = U^{\#}_i/G \subset \mathbb{P}$.
Then we have a factorization of $\pi'_i = \pi' \vert_{U^{\#}_i}$:
\[
\begin{CD}
\pi'_i: U^{\#}_i @>{\pi_i}>> \mathcal{U}_i @>{\sigma_i}>> U_i = 
\mathcal{U}_i/\mu_{a_i}
\end{CD}
\]
where $\mathcal{U}_i$ is smooth and $\sigma_i$ is etale in codimension $1$.
Thus the set of coverings $\{\sigma_i\}$ define the stack $\mathcal{P}$.

For $i \ne j$, if we set $U^{\#}_{ij} = U^{\#}_i \cap U^{\#}_j$ 
and $U_{ij} = U^{\#}_{ij}/G \subset \mathbb{P}$,
then $\pi'_{ij} = \pi' \vert_{U^{\#}_{ij}}$ is factorized as 
\[
\begin{CD}
\pi'_{ij}: U^{\#}_{ij} @>{\pi_{ij}}>> \mathcal{U}_{ij} @>{\sigma_{ij}}>> 
U_{ij} = \mathcal{U}_{ij}/(\mu_{a_i} \times \mu_{a_j})
\end{CD}
\]
where $\mathcal{U}_{ij} = \mathcal{U}_i \times_{\mathcal{P}} \mathcal{U}_j$.

We have 
\[
\pi^*\mathcal{O}_{\mathcal{P}}(1) \cong 
\mathcal{O}_{\mathbb{P}^{\#}}(1).
\]
The images $\mathcal{D}_i = \pi(D^{\#}_i)$ are prime divisors on 
$\mathcal{P}$ such that
$\pi^*\mathcal{D}_i = a_iD^{\#}_i$ and
\[
\mathcal{O}_{\mathbb{P}^{\#}}(D^{\#}_i) 
\cong \mathcal{O}_{\mathbb{P}^{\#}}(1), \quad 
\mathcal{O}_{\mathcal{P}}(\mathcal{D}_i) \cong \mathcal{O}_{\mathcal{P}}(a_i).
\] 
The dualizing invertible sheaf is given by 
\[
\omega_{\mathcal{P}} \cong \mathcal{O}_{\mathcal{P}}(\sum_i a_i).
\]

We have the following resolution of the structure sheaf of the diagonal in 
$\mathbb{P}^{\#} \times \mathbb{P}^{\#}$ due to Beilinson 
\cite{Beilinson}.
From a tautological exact sequence on $\mathbb{P}^{\#}$
\[
0 \to \mathcal{O}_{\mathbb{P}^{\#}}(-1) \to 
\mathcal{O}_{\mathbb{P}^{\#}}^{n+1} \to T_{\mathbb{P}^{\#}}(-1) \to 0
\]
we obtain a homomorphism of sheaves 
on $\mathbb{P}^{\#} \times \mathbb{P}^{\#}$
\[
s: \mathcal{O}_{\mathbb{P}^{\#}}(-1) \boxtimes 
\mathcal{O}_{\mathbb{P}^{\#}}
\to \mathcal{O}_{\mathbb{P}^{\#}} \boxtimes
T_{\mathbb{P}^{\#}}(-1)
\]
given by
\[
s(P,Q)(v) = v \text{ mod } \mathbb{C}w 
\text{ if } [v] = P \text{ and } [w] = Q
\]
where $v,w \in V$ for $\mathbb{P}^{\#} = \mathbb{P}(V^*)$.
We can also write
\[
s = \sum_{i=0}^n x^{\#}_i \boxtimes \frac{\partial}{\partial y^{\#}_i}.
\]
The cokernel of the homomorphism
\[
s': \mathcal{O}_{\mathbb{P}^{\#}}(-1) \boxtimes 
\Omega^1_{\mathbb{P}^{\#}}(1)
\to \mathcal{O}_{\mathbb{P}^{\#}} \boxtimes
\mathcal{O}_{\mathbb{P}^{\#}}
\]
induced by $s$ is the structure sheaf of the diagonal
$\Delta \mathbb{P}^{\#} \subset 
\mathbb{P}^{\#} \times \mathbb{P}^{\#}$.
We note that $s'$ is equivariant under the action of the diagonal subgroup 
$\Delta G \subset G \times G$.

Therefore, we obtain an exact sequence called {\it Beilinson resolution} 
of the diagonal as a Koszul complex
\begin{equation}\label{Beilinson}
\begin{split}
&0 \to \mathcal{O}_{\mathbb{P}^{\#}}(-n) \boxtimes 
\Omega^n_{\mathbb{P}^{\#}}(n) \to
\mathcal{O}_{\mathbb{P}^{\#}}(-n+1) \boxtimes 
\Omega^{n-1}_{\mathbb{P}^{\#}}(n-1) \to \\
&\dots \to \mathcal{O}_{\mathbb{P}^{\#}}(-1) \boxtimes 
\Omega^1_{\mathbb{P}^{\#}}(1) \to 
\mathcal{O}_{\mathbb{P}^{\#}} \boxtimes 
\mathcal{O}_{\mathbb{P}^{\#}} \to 
\mathcal{O}_{\Delta \mathbb{P}^{\#}} \to 0.
\end{split}
\end{equation}
We note that the homomorphisms are $\Delta G$-equivariant.

We extend the Beilinson resolution for the stacky sheaves on
weighted projective spaces according to the description 
by Canonaco \cite{Canonaco}, where he considered the sheaves on $\mathbb{P}$ 
instead of $\mathcal{P}$.

The group of characters $G^*$ of $G$ is isomorphic to 
$\mathbb{Z}_{a_0} \times \dots \mathbb{Z}_{a_n}$.
For a character $\chi = (\chi_0, \dots, \chi_n) \in G^*$, we take 
representatives of the $\chi_i$ such that $0 \le \chi_i < a_i$, 
and define $\vert \chi \vert = \sum_{i=0}^n \chi_i \in \mathbb{Z}$.

Since $\pi$ is ramified along the $\mathcal{D}_i$ with ramification 
indices $a_i$,
we have a decomposition of the direct image sheaf 
into eigenspaces according to the $G$-action:
\[
\pi_*\mathcal{O}_{\mathbb{P}^{\#}}
\cong \bigoplus_{\chi} (\pi_*\mathcal{O}_{\mathbb{P}^{\#}})^{\chi}, \quad 
(\pi_*\mathcal{O}_{\mathbb{P}^{\#}})^{\chi} 
\cong \mathcal{O}_{\mathcal{P}}(- \vert \chi \vert).
\]
By \cite{Canonaco}, we have also 
\[
\pi_*\Omega^p_{\mathbb{P}^{\#}}
\cong \bigoplus_{\chi} (\pi_*\Omega^p_{\mathbb{P}^{\#}})^{\chi}, \quad  
(\pi_*\Omega^p_{\mathbb{P}^{\#}})^{\chi} 
\cong \Omega^p_{\mathcal{P}}(\log \mathcal{D}_{\chi})(- \vert \chi \vert)
\]
where $\mathcal{D}_{\chi} = \sum_{\chi_i \ne 0} \mathcal{D}_i$.

Let $a$ be a coherent sheaf on $\mathcal{P}$, and consider the 
Beilinson resolution of the pull-back $a^{\#} = \pi^*a$, where we 
note that $\pi$ is a flat morphism because we consider
the stack $\mathcal{P}$ instead of the variety $\mathbb{P}$.
We have a quasi-isomorphism 
\begin{equation}\label{resol}
\begin{split}
&\mathcal{O}_{\Delta \mathbb{P}^{\#}} \otimes p_2^*a^{\#} 
\cong \{0 \to \mathcal{O}_{\mathbb{P}^{\#}}(-n) \boxtimes 
(\Omega^n_{\mathbb{P}^{\#}}(n) \otimes a^{\#}) \to \\
&\dots \to \mathcal{O}_{\mathbb{P}^{\#}}(-1) \boxtimes 
(\Omega^1_{\mathbb{P}^{\#}}(1) \otimes a^{\#}) \to 
\mathcal{O}_{\mathbb{P}^{\#}} \boxtimes 
a^{\#} \to 0\}.
\end{split}
\end{equation}
Since 
\[
p_{1*}(\mathcal{O}_{\Delta \mathbb{P}^{\#}} \otimes p_2^*a^{\#})
\cong a^{\#}
\]
we have
\[
\begin{split}
&a^{\#} \cong 
\text{sg}\{0 \to \mathcal{O}_{\mathbb{P}^{\#}}(-n) \otimes_{\mathbb{C}} 
R\Gamma(\mathbb{P}^{\#}, \Omega^n_{\mathbb{P}^{\#}}(n) 
\otimes a^{\#}) \to \\
&\dots \to \mathcal{O}_{\mathbb{P}^{\#}}(-1) \otimes_{\mathbb{C}} 
R\Gamma(\mathbb{P}^{\#}, \Omega^1_{\mathbb{P}^{\#}}(1) \otimes a^{\#}) 
\to \mathcal{O}_{\mathbb{P}^{\#}} \otimes_{\mathbb{C}} 
R\Gamma(\mathbb{P}^{\#}, a^{\#}) \to 0\}
\end{split}
\]
where $\text{sg}$ stands for the associated single complex and 
\[
R\Gamma(\mathbb{P}^{\#}, \Omega^p_{\mathbb{P}^{\#}}(p) \otimes a^{\#})
= \sum_q H^q(\mathbb{P}^{\#}, \Omega^p_{\mathbb{P}^{\#}}(p) 
\otimes a^{\#})[-q].
\]
Since the group $G \cong \Delta G$ acts equivariantly on $(\ref{resol})$,
we have
\[
\begin{split}
&a \cong 
\text{sg}\{0 \to \bigoplus_{\chi} \mathcal{O}_{\mathcal{P}}
(- n - \vert \chi \vert) \otimes 
R\Gamma(\mathcal{P}, \Omega^n_{\mathcal{P}}(\log \mathcal{D}_{- \chi})
(n - \vert - \chi \vert) \otimes a) \to \\
&\dots \to \bigoplus_{\chi} \mathcal{O}_{\mathcal{P}}
(- 1 - \vert \chi \vert) \otimes 
R\Gamma(\mathcal{P}, \Omega^1_{\mathcal{P}}(\log \mathcal{D}_{- \chi})
(1 - \vert - \chi \vert) \otimes a) \\
&\to \bigoplus_{\chi} \mathcal{O}_{\mathcal{P}}(- \vert \chi \vert)
\otimes R\Gamma(\mathcal{P}, a(- \vert - \chi \vert)) \to 0\}.
\end{split}
\]

In particular, we have a spectral sequence of Beilinson type
\[
E_1^{p,q} = \bigoplus_{\chi} \mathcal{O}_{\mathcal{P}}(p - \vert \chi \vert)
\otimes H^q(\mathcal{P}, \Omega^{-p}_{\mathcal{P}}
(\log \mathcal{D}_{- \chi})(- p - \vert - \chi \vert) \otimes a)
\Rightarrow a
\]
where the right hand side is at degree $0$.

If $H^q(\mathcal{P}, \Omega^p_{\mathcal{P}}(\log \mathcal{D}_{- \chi})
(p - \vert - \chi \vert) \otimes a) = 0$
for $q > 0$ and all $p$ and $\chi$, 
then we obtain an exact sequence
\[
\begin{split}
&0 \to \bigoplus_{\chi} \mathcal{O}_{\mathcal{P}}(- n - \vert \chi \vert) 
\otimes H^0(\mathcal{P}, \Omega^n_{\mathcal{P}}(\log \mathcal{D}_{- \chi})
(n - \vert - \chi \vert) \otimes a) \to 
\\
&\dots \to \bigoplus_{\chi} \mathcal{O}_{\mathcal{P}}(- 1 - \vert \chi \vert) 
\otimes H^0(\mathcal{P}, \Omega^1_{\mathcal{P}}(\log \mathcal{D}_{- \chi})
(1 - \vert - \chi \vert) \otimes a) \\
&\to \bigoplus_{\chi} \mathcal{O}_{\mathcal{P}}(- \vert \chi \vert) \otimes 
H^0(\mathcal{P}, a(- \vert - \chi \vert)) \to a \to 0.
\end{split}
\]
We call this the {\it left resolution} of $a$, because the direction of the 
resolution is to the left.

On the other hand, 
if $H^q(\mathcal{P}, \Omega^p_{\mathcal{P}}(\log \mathcal{D}_{- \chi})
(p - \vert - \chi \vert) \otimes a) = 0$
for $q < n$ and all $p$ and $\chi$, then we obtain an exact sequence
called the {\it right resolution} of $a$:
\[
\begin{split}
&0 \to a \to 
\bigoplus_{\chi} \mathcal{O}_{\mathcal{P}}(- n - \vert \chi \vert) \otimes 
H^n(\mathcal{P}, \Omega^n_{\mathcal{P}}(\log \mathcal{D}_{- \chi})
(n - \vert - \chi \vert) \otimes a) \to 
\\
&\dots \to \bigoplus_{\chi} \mathcal{O}_{\mathcal{P}}(- 1 - \vert \chi \vert) 
\otimes H^n(\mathcal{P}, \Omega^1_{\mathcal{P}}(\log \mathcal{D}_{- \chi})
(1 - \vert - \chi \vert) \otimes a) \\
&\to \bigoplus_{\chi} \mathcal{O}_{\mathcal{P}}(- \vert \chi \vert) \otimes 
H^n(\mathcal{P}, a(- \vert - \chi \vert)) \to 0.
\end{split}
\]
We note that the left and right resolutions are related by the shift.

If we define an object $e \in D^b(\text{Coh}(\mathcal{P} \times \mathcal{P}))$ 
by 
\[
\begin{split}
&e = 
\{0 \to \bigoplus_{\chi} \mathcal{O}_{\mathcal{P}}
(- n - \vert \chi \vert) \boxtimes 
\Omega^n_{\mathcal{P}}(\log \mathcal{D}_{- \chi})(n - \vert - \chi \vert) \to 
\\
&\dots \to \bigoplus_{\chi} \mathcal{O}_{\mathcal{P}}
(- 1 - \vert \chi \vert) \boxtimes 
\Omega^1_{\mathcal{P}}(\log \mathcal{D}_{- \chi})(1 - \vert - \chi \vert)  \\
&\to \bigoplus_{\chi} \mathcal{O}_{\mathcal{P}}(- \vert \chi \vert)
\boxtimes \mathcal{O}_{\mathcal{P}}(- \vert - \chi \vert) \to 0\}
\end{split}
\]
then we have
\[
a \cong p_{1*}(e \otimes p_2^*a) \cong p_{2*}(e \otimes p_1^*a).
\]


\section{Proof of Theorem~\ref{main}}

The proof of the theorem is partly parallel to the original proof in 
\cite{Orlov}.  But we cannot use an embedding in a projective space 
because $\mathcal{X}$ has more sheaves than $X$.

We denote by $\sigma: \mathcal{X} \to X$ the natural morphism, and
$\pi: X^{\#} \to \mathcal{X}$ the covering considered in the last section.
We fix a very ample line bundle $L$ on $X$ such that the pull-back 
$\pi^*\mathcal{L}$ for $\mathcal{L} = \sigma^*L$
on the covering $X^{\#}$ has the property described in 
Theorem~\ref{diagonal1}.

\subsection{Step 1}

We prove the boundedness of the functor $F$ (\cite{Orlov}~Lemma~2.4), that is,
there exists a fixed interval of integers such that 
$H^k(F(\mathcal{A})) = 0$ for any integer $k$ outside this interval and 
for any coherent sheaf $\mathcal{A}$ on $\mathcal{X}$.

Let $F^*$ be the left adjoint functor of $F$.
Fix an embedding $Y \to \mathbb{P}^M$ into some projective space.
By pulling back the Beilinson resolution to $\mathcal{Y}$, 
we obtain right resolutions
\[
\mathcal{O}_{\mathcal{Y}}(- j) \cong
\{V_M^j \otimes_{\mathbb{C}} \mathcal{O}_{\mathcal{Y}}
\to V_{M-1}^j \otimes_{\mathbb{C}} \mathcal{O}_{\mathcal{Y}}(1)
\to \dots \to V_0^j \otimes_{\mathbb{C}} \mathcal{O}_{\mathcal{Y}}(M)\}
\]
for integers $j > 0$, where 
$V_p^j = H^M(\mathbb{P}^M, \Omega^p(p - j - M))$ ($0 \le p \le M$).

Let $\mathcal{B}_0$ be a locally free sheaf on $\mathcal{Y}$
obtained in Theorem~\ref{locallyfree}.
We choose integers $k_1 < k_2$ such that 
\[
H^k(F^*(\mathcal{B}_0(j))) = 0 
\]
for $k \not\in [k_1,k_2]$ and $0 \le j \le M$.
Then 
\[
H^k(F^*(\mathcal{B}_0(- j))) = 0 
\]
for $k \not\in [k_1,k_2+M]$ and any $j > 0$.

For any coherent sheaf $\mathcal{A}$ on $\mathcal{X}$, we have 
\[
\text{Hom}^k(\mathcal{B}_0(- j), F(\mathcal{A})) 
\cong \text{Hom}^k(F^*(\mathcal{B}_0(- j)), \mathcal{A}) \cong 0
\]
for $k \not\in [-k_2-M, -k_1+\dim X]$ and for any $j > 0$.
If we take $j$ sufficiently large, then we have
\[
\text{Hom}^p(\mathcal{B}_0(- j), H^q(F(\mathcal{A}))) \cong 0
\]
for $p > 0$ and any $q$.
Hence $H^k(F(\mathcal{A})) = 0$ for $k \not\in [-k_2-M, -k_1+\dim X]$.

By replacing $F$ by its shift, we assume from now on that 
there exists an integer $k_0$ such that 
$H^k(F(\mathcal{A})) = 0$ for $k \not\in [-k_0, 0]$ and any sheaf 
$\mathcal{A}$ on $\mathcal{X}$.


\subsection{Step 2}

We shall define an object 
$e \in D^b(\text{Coh}(\mathcal{X} \times \mathcal{Y}))$.

We fix an integer $m$ such that $m > \dim X + \dim Y + k_0$, 
where the length $k_0$ is given in Step 1, and let 
$c_{\bullet}$ be the complex defined by
\begin{equation}\label{e(m)}
c_p = \begin{cases}
(\pi_*\pi^*\mathcal{L}^{-p} \boxtimes F(\pi_*R_p))^{\Delta G} 
&\text{ for } 0 \le p \le m \\
0 &\text{ otherwise}
\end{cases}
\end{equation}
with the morphisms $\delta_p: c_p \to c_{p-1}$ induced from $(\ref{LRhom})$.

We have
\[
\begin{split}
&\text{Hom}(\pi_*\pi^*\mathcal{L}^{-p} \boxtimes F(\pi_*R_p)[r],
\pi_*\pi^*\mathcal{L}^{-q} \boxtimes F(\pi_*R_q)) \\
&\cong \bigoplus_{r_1+r_2=r}
\text{Hom}(\pi_*\pi^*\mathcal{L}^{-p}[r_1], \pi_*\pi^*\mathcal{L}^{-q})
\otimes \text{Hom}(F(\pi_*R_p)[r_2], F(\pi_*R_q)) \\
&\cong \bigoplus_{r_1+r_2=r}
\text{Hom}(\pi_*\pi^*\mathcal{L}^{-p}[r_1], \pi_*\pi^*\mathcal{L}^{-q})
\otimes \text{Hom}(\pi_*R_p[r_2], \pi_*R_q) \\
&\cong 0
\end{split}
\]
if $r > 0$.
Thus there exists a right convolution 
$\delta_0: c_0 \to e' \in D^b(\text{Coh}(\mathcal{X} \times \mathcal{Y}))$.

\begin{Lem}\label{e-split}
$H^p(e') = 0$ unless $p \in [-m-k_0, -m] \cup [-k_0, 0]$.
\end{Lem}

\begin{proof}
Assume that there exists $p_0 \not\in [-m-k_0, -m] \cup [-k_0, 0]$ such that 
$H^{p_0}(e') \ne 0$.
For a locally free sheaf $\mathcal{A}$ on $\mathcal{X}$, 
we have a spectral sequence
\[
E_2^{p,q} = R^pp_{2*}(H^q(e') \otimes p_1^*\mathcal{A}) \Rightarrow
H^{p+q}(p_{2*}(e' \otimes p_1^*\mathcal{A})).
\]
Hence there exists $\mathcal{A}$ such that
$H^{p_0}(p_{2*}(e' \otimes p_1^*\mathcal{A})) \ne 0$ by 
Theorem~\ref{locallyfree}. 
We may also assume that 
\begin{equation}\label{condA}
H^p(\mathcal{X}, \mathcal{A} \otimes \pi_*\pi^*\mathcal{L}^{-q}) = 0
\text{ for } p > 0, 0 \le q \le m + \dim X.
\end{equation}

We consider a complex of sheaves $a_{\bullet}$ given by
\[
a_p = \begin{cases} 
(H^0(\mathcal{X}, \mathcal{A} \otimes \pi_*\pi^*\mathcal{L}^{-p}) 
\otimes \pi_*R_p)^{\Delta G} 
&\text{ for } 0 \le p \le m \\
0 &\text{ otherwise.} 
\end{cases}
\]
This coincides with $p_{2*}(e^- \otimes p_1^*\mathcal{A})$ up to degree $-m$
for $e^-$ given in Theorem~\ref{diagonal2}.

If we denote $U(\mathcal{A}) = \text{Ker}(a_m \to a_{m-1})$, then
there is a distinguished triangle
\[
U(\mathcal{A})[m] \to a_{\bullet} \to \mathcal{A} \to U(\mathcal{A})[m+1].
\]
Since $\text{Hom}_{D^b(\mathcal{X})}(\mathcal{A}, U(\mathcal{A})[m+1]) = 0$, 
we have a right convolution
$d_0: a_0 \to U(\mathcal{A})[m] \oplus \mathcal{A}$ of $a_{\bullet}$.
Hence $F(d_0): F(a_0) \to F(U(\mathcal{A})[m] \oplus \mathcal{A})$ is also 
a right convolution of the complex of objects $F(a_{\bullet})$.

On the other hand, since we have 
\[
p_{2*}(c_p \otimes p_1^*\mathcal{A}) = F(a_p)
\]
$\delta_0$ induces a right convolution
$\delta_{0*}: F(a_0) \to p_{2*}(e' \otimes p_1^*\mathcal{A})$
of the complex $F(a_{\bullet})$.

Since both $F(U(\mathcal{A})[m] \oplus \mathcal{A})$ and 
$p_{2*}(e' \otimes p_1^*\mathcal{A})$
are right convolutions of the same complex, there exists an isomorphism
$\tilde f(\mathcal{A})$, which is not uniquely determined, 
making the following diagram commutative
\[
\begin{CD}
F(a_0) @>{\bar d_{0*}}>> p_{2*}(e' \otimes p_1^*\mathcal{A}) \\
@V=VV @V{\tilde f(\mathcal{A})}VV \\
F(a_0) @>{F(d_0)}>> F(U(\mathcal{A})[m] \oplus \mathcal{A}).
\end{CD}
\]
Therefore, we obtain a contradiction by Step 1.
\end{proof}

It follows that there exist objects $e, e'' \in 
D^b(\text{Coh}(\mathcal{X} \times \mathcal{Y}))$ and a distinguished triangle
\[
e'' \to e' \to e \to e''[1]
\]
such that $H^p(e) = 0$ (resp. $H^p(e'') = 0$) unless
$p \in [-k_0, 0]$ (resp. $[-m-k_0, -m]$).
Since $\text{Hom}_{D^b(\text{Coh}(\mathcal{X} \times \mathcal{Y}))}
(e, e''[1]) = 0$, 
we have $e' \cong e \oplus e''$.


\subsection{Step 3}

We construct a functorial isomorphism 
from $\Phi^e_{\mathcal{X} \to \mathcal{Y}}$ to $F$ 
for a certain set of locally free sheaves.

\begin{Lem}\label{functorial1}
Let $\mathcal{A}$ be any locally free coherent sheaf on $\mathcal{X}$ 
such that 
\[
H^p(\mathcal{X}, \mathcal{A} \otimes \pi_*\pi^*\mathcal{L}^{-q}) = 0
\] 
for any $p > 0$ and
$0 \le q \le m + \dim X$.
Then there exists an isomorphism
\[
f(\mathcal{A}): \Phi^e_{\mathcal{X} \to \mathcal{Y}}(\mathcal{A}) 
\to F(\mathcal{A})
\] 
which is functorial in the sense that,
for any homomorphism $\alpha: \mathcal{A} \to \mathcal{B}$ 
to another locally free coherent sheaf $\mathcal{B}$
which satisfies the same condition, the following diagram is commutative
\[
\begin{CD}
\Phi^e_{\mathcal{X} \to \mathcal{Y}}(\mathcal{A}) @>{f(\mathcal{A})}>> 
F(\mathcal{A}) \\
@V{\Phi^e_{\mathcal{X} \to \mathcal{Y}}(\alpha)}VV  @V{F(\alpha)}VV \\
\Phi^e_{\mathcal{X} \to \mathcal{Y}}(\mathcal{B}) @>{f(\mathcal{B})}>> 
F(\mathcal{B}).
\end{CD}
\]
\end{Lem}

\begin{proof}
We define $a_{\bullet}$, $U(\mathcal{A})$, $d_0: a_0 \to U(\mathcal{A})[m] 
\oplus \mathcal{A}$, 
$\delta_{0*}: F(a_0) \to p_{2*}(e' \otimes p_1^*\mathcal{A})
= \Phi^{e'}_{\mathcal{X} \to \mathcal{Y}}(\mathcal{A})$ and 
$\tilde f(A)$ as in Lemma~\ref{e-split}.

The isomorphism $\tilde f(\mathcal{A})$ induces an isomorphism 
$p_{2*}(e \otimes p_1^*\mathcal{A}) \to F(\mathcal{A})$ 
by Steps 1 and 2.  But we note that it is not sufficient, because 
$\tilde f(\mathcal{A})$ is not uniquely determined.

Let $\epsilon(\mathcal{A}): F(U(\mathcal{A})[m] \oplus \mathcal{A}) 
\to F(\mathcal{A})$ be the projection.
Since 
\[
\text{Hom}(F(a_p)[r], F(\mathcal{A})) \cong 
\text{Hom}(a_p[r], \mathcal{A}) \cong 0
\]
for any $p$ and $r > 0$,
$\epsilon(\mathcal{A})$ is the only morphism which makes 
the following diagram commutative by Lemma~\ref{complexhom1lem}:
\[
\begin{CD}
F(a_0) @>{F(d_0)}>> F(U(\mathcal{A})[m] \oplus \mathcal{A}) \\
@V=VV @V{\epsilon(\mathcal{A})}VV \\
F(a_0) @>{\epsilon(\mathcal{A}) F(d_0)}>> F(\mathcal{A}).
\end{CD}
\]

By the same reason as above, the projection
$\epsilon'(\mathcal{A}): p_{2*}(e' \otimes p_1^*\mathcal{A})
\to p_{2*}(e \otimes p_1^*\mathcal{A})$ is the only morphism
making the following diagram commutative
\[
\begin{CD}
F(a_0) @>{\delta_{0*}}>> p_{2*}(e' \otimes p_1^*\mathcal{A}) \\
@V=VV @V{\epsilon'(\mathcal{A})}VV \\
F(a_0) @>{\epsilon'(\mathcal{A})\delta_{0*}}>> 
p_{2*}(e \otimes p_1^*\mathcal{A}).
\end{CD}
\]
Therefore, there exists a uniquely determined isomorphism 
$f(\mathcal{A}): p_{2*}(e \otimes p_1^*\mathcal{A}) \to F(\mathcal{A})$ 
such that the following diagram is commutative
\[
\begin{CD}
F(a_0) @>{\epsilon'(\mathcal{A})\delta_{0*}}>> 
p_{2*}(e \otimes p_1^*\mathcal{A}) \\
@V=VV @V{f(\mathcal{A})}VV \\
F(a_0) @>{\epsilon(\mathcal{A}) F(d_0)}>> F(\mathcal{A}).
\end{CD}
\]

Next, we consider a complex of sheaves $b_{\bullet}$ given by 
\[
b_p = \begin{cases} 
(H^0(\mathcal{X}, \mathcal{B} \otimes \pi_*\pi^*\mathcal{L}^{-p}) 
\otimes \pi_*R_p)^{\Delta G} 
&\text{ for } 0 \le p \le m \\
0 &\text{ otherwise}. 
\end{cases}
\]
The homomorphism $\alpha: \mathcal{A} \to \mathcal{B}$ induces a morphism of 
complexes $\alpha_{\bullet}: a_{\bullet} \to b_{\bullet}$.
Let $e_0: b_0 \to U(\mathcal{B})[m] \oplus \mathcal{B}$ be the 
right convolution. 

There exists a uniquely determined morphism $g$
to make the following diagram commutative
\[
\begin{CD}
F(a_0) @>{\delta_{0*}}>> p_{2*}(e' \otimes p_1^*\mathcal{A}) \\
@V{F(\alpha_0)}VV @V{g}VV \\
F(b_0) @>{\epsilon(\mathcal{B}) F(e_0)}>> F(\mathcal{B}).
\end{CD}
\]
It follows that 
\[
g = f(\mathcal{B})\Phi^e_{\mathcal{X} \to \mathcal{Y}}(\alpha)
\epsilon'(\mathcal{A})
= F(\alpha) f(\mathcal{A}) \epsilon'(\mathcal{A})
\]
hence $f(\mathcal{B}) \Phi^e_{\mathcal{X} \to \mathcal{Y}}(\alpha)
= F(\alpha) f(\mathcal{A})$.
\end{proof}


\subsection{Step 4}

We extend the functorial isomorphism from $F$ to 
$\Phi^e_{\mathcal{X} \to \mathcal{Y}}$ to an ample set of objects.

A sequence of objects $\{P_k\}_{k \in \mathbb{Z}}$ in an abelian 
category $\mathcal{C}$ is said to be {\it ample}
if for any object $C \in \mathcal{C}$, there exists an integer $k_0(C)$
such that the following conditions are satisfied for any $k < k_0(C)$
(\cite{Orlov}~Definition~2.12).

(1) The canonical morphism $\text{Hom}(P_k, C) \otimes P_k \to C$ is 
surjective.

(2) $\text{Hom}(P_k, C[r]) = 0$ if $r \ne 0$.

(3) $\text{Hom}(C, P_k) = 0$.

We note that the set $\{P_k\}_{k \in \mathbb{Z}}$ 
becomes a spanning class 
of the derived category $D^b(\mathcal{C})$ if it has a Serre functor
(cf. \cite{Bridgeland}).

By Theorem~\ref{locallyfree}, we obtain the following:

\begin{Lem}
There exists a locally free sheaf $\mathcal{A}_0$ on $\mathcal{X}$ such that 
the sequence $\{\mathcal{A}_0 \otimes \mathcal{L}^k\}$ 
is ample in the category 
of coherent sheaves $\text{Coh}(\mathcal{X})$.
\end{Lem}

\begin{Lem}\label{functorial2}
There exist isomorphisms 
$f_k: F(\mathcal{A}_0 \otimes \mathcal{L}^k) \to 
\Phi^e_{\mathcal{X} \to \mathcal{Y}}(\mathcal{A}_0 \otimes \mathcal{L}^k)$ 
for any integers $k$ which are functorial in the sense that,
for any homomorphism $\alpha: \mathcal{A}_0 \otimes \mathcal{L}^k \to 
\mathcal{A}_0 \otimes \mathcal{L}^{k'}$, 
the following diagram is commutative
\[
\begin{CD}
F(\mathcal{A}_0 \otimes \mathcal{L}^k) @>{f_k}>> 
\Phi^e_{\mathcal{X} \to \mathcal{Y}}(\mathcal{A}_0 \otimes \mathcal{L}^k) \\
@V{F(\alpha)}VV  @V{\Phi^e_{\mathcal{X} \to \mathcal{Y}}(\alpha)}VV \\
F(\mathcal{A}_0 \otimes \mathcal{L}^{k'}) @>{f_{k'}}>> 
\Phi^e_{\mathcal{X} \to \mathcal{Y}}(\mathcal{A}_0 \otimes \mathcal{L}^{k'}).
\end{CD}
\]
\end{Lem}

\begin{proof}
We denote $\Phi = \Phi^e_{\mathcal{X} \to \mathcal{Y}}$.
By Lemma~\ref{functorial1}, there exists an integer $k_1$ such that our
assertion holds if $k, k' \ge k_1$.
We proceed by the descending induction on such $k_1$. 

Let us fix an embedding $X \to \mathbb{P}^N$ such that 
$L = \mathcal{O}_X(1)$.
We have an exact sequence
\begin{equation}\label{resolutionL^-1}
0 \to \mathcal{O}_{\mathcal{X}} \to V_N \otimes \mathcal{L} \to 
V_{N-1} \otimes \mathcal{L}^2 \to \dots \to V_0 \otimes \mathcal{L}^{N+1} \to 0
\end{equation}
for $V_p = H^N(\mathbb{P}^N, \Omega^p_{\mathbb{P}^N}(p-N-1))$.

For any integer $k$, we consider a complex of objects $a_{\bullet}^{(k)}$
and $b_{\bullet}^{(k)}$ defined by 
\[
\begin{split}
&a_p^{(k)} = \begin{cases} F(V_p \otimes \mathcal{A}_0 \otimes 
\mathcal{L}^{k+N+1-p}) 
&\text{ for } 0 \le p \le N \\
0 &\text{ otherwise} 
\end{cases} \\
&b_p^{(k)} = \begin{cases} \Phi(V_p \otimes \mathcal{A}_0 \otimes 
\mathcal{L}^{k+N+1-p}) 
&\text{ for } 0 \le p \le N \\
0 &\text{ otherwise.} 
\end{cases}
\end{split}
\]
If $k \ge k_1 -1$, then there is an isomorphism of complexes 
$f^{(k)}_{\bullet}: a_{\bullet}^{(k)} \to b_{\bullet}^{(k)}$ 
by the induction hypothesis, where 
\[
f^{(k)}_p = \text{Id}_{V_p} \otimes f_{k+N+1-p}.
\]
We have the left convolutions 
$d^{(k)}: F(\mathcal{A}_0 \otimes \mathcal{L}^k) \to a_N^{(k)}$ 
and 
$e^{(k)}: \Phi(\mathcal{A}_0 \otimes \mathcal{L}^k) \to b_N^{(k)}$, 
and there is a uniquely determined isomorphism
$f_k: F(\mathcal{A}_0 \otimes \mathcal{L}^k) \to 
\Phi(\mathcal{A}_0 \otimes \mathcal{L}^k)$ 
such that 
$e^{(k)}f_k = f^{(k)}_Nd^{(k)}$.
Note that if $k \ge k_1$, then the morphism $f_k$ which we already have 
satisfies this condition thanks to the functoriality of the $f_k$.

Let $k, k' \ge k_1 - 1$.  
Then the homomorphism $\alpha$ induces a homomorphism of complexes
$\alpha_{*\bullet}: a_{\bullet}^{(k)} \to b_{\bullet}^{(k')}$
given by 
\[
\alpha_{*p} = \text{Id}_{V_p} \otimes f_{k'+N+1-p}
F(\alpha \otimes \text{Id}_{\mathcal{L}^{N+1-p}}).
\]
Therefore, we have a uniquely determined homomorphism
$g: F(\mathcal{A}_0 \otimes \mathcal{L}^k) \to 
\Phi(\mathcal{A}_0 \otimes \mathcal{L}^{k'})$ 
by the condition that 
$e^{(k')}g = \alpha_{*N}d^{(k)}$.
Since both $\Phi(\alpha)f_k$ and $f_{k'}F(\alpha)$
satisfiy this condition, they coincide.
\end{proof}


\subsection{Step 5}

By \cite{Orlov}~Proposition~2.16, we obtain an isomorphism of 
functors $f: F \to \Phi^e_{\mathcal{X} \to \mathcal{Y}}$ 
by extending the $f_k$: 

\begin{Lem}
Let $\mathcal{C}$ be an abelian category, $D^b(\mathcal{C})$ its derived 
category, and $\mathcal{D}$ another triangulated category.
Let $F: D^b(\mathcal{C}) \to \mathcal{D}$ be a fully faithful exact functor
which has a left adjoint, and $G: D^b(\mathcal{C}) \to \mathcal{D}$ 
another exact functor which has a left adjoint.
Assume that $D^b(\mathcal{C})$ and $\mathcal{D}$ have Serre functors and that
there is an ample sequence $\Omega = \{P_k\}$ in $\mathcal{C}$
with an isomorphism of functors $f_{\Omega}: F \vert_{\Omega}
\to G \vert_{\Omega}$ when restricted to the full subcategory 
$\Omega \subset D^b(\mathcal{C})$.
Then there exists an isomorphism of functors $f: F \to G$
which is an extension of $f_{\Omega}$. 
\end{Lem}

\begin{proof}
Let $F^*$ and $F^!$ (resp. $G^*$ and $G^!$) be the left and right 
adjoint functors of $F$ (resp. $G$).
Since $F$ is fully faithful, the natural morphism of functors
$\text{Id}_{D^b(\mathcal{C})} \to F^!F$ is an isomorphism,
because 
\[
\text{Hom}(a, b) \to \text{Hom}(a, F^!Fb) 
\]
is an isomorphism for any objects $a$ and $b$.

Since $G$ is isomorphic to $F$ when restricted to a spanning class 
$\Omega$, it is also fully faithful, 
hence the natural morphism of functors 
$G^*G \to \text{Id}_{D^b(\mathcal{C})}$ is an isomorphism.

$G^* F$ is the left adjoint functor of $F^!G$.
Since they are isomorphic to the identity functor when restricted to 
$\Omega$, it follows that they are fully faithful, 
hence they are quasi-inverses of each other.
Thus $F^!G$ is an autoequivalence of $D^b(\mathcal{C})$, and 
there exists an extended isomorphism of functors
$f': \text{Id}_{D^b(\mathcal{C})} \to F^!G$ by 
\cite{Orlov}~Proposition~2.16. 
By adjunction, we obtain a morphism of functors
$f: F \to G$.

We prove that $f$ is an isomorphism.
For any object $a$, let 
\[
\begin{CD}
F(a) @>{f(a)}>> G(a) @>>> c @>>> F(a)[1] 
\end{CD}
\]
be a distinguished triangle.
Since $F^!(f(a))$ is an isomorphism, we have $F^!(c) \cong 0$.
From
\[
\begin{split}
&\text{Hom}(\omega, G^!(c)) \cong \text{Hom}(G(\omega), c) \\
&\cong \text{Hom}(F(\omega), c) \cong \text{Hom}(\omega, F^!(c)) \cong 0
\end{split}
\]
for any $\omega \in \Omega$, we deduce that $G^!(c) \cong 0$.
Hence $\text{Hom}(G(a), c) \cong 0$, 
thus $F(a) \cong G(a) \oplus c[-1]$.
But 
\[
\text{Hom}(F(a), c[-1]) \cong \text{Hom}(a, F^!(c)[-1]) \cong 0
\]
hence $c = 0$.
\end{proof}


\subsection{Step 6}

Let $e_1 \in D^b(\text{Coh}(\mathcal{X} \times \mathcal{Y}))$ 
be any object such that 
there is an isomorphism of functors 
$F \to \Phi^{e_1}_{\mathcal{X} \to \mathcal{Y}}$.
We shall prove that $e_1 \cong e$.

Let $(c_{\bullet}, \delta_{\bullet})$ be the complex of objects 
considered in Step 2.
We consider a commutative diagram
\[
\begin{CD}
c_1 @>>> c_{10} @>>> c_0 \\
@VVV @VVV @VVV \\
\bar c_1 @>>> \bar c_{10} @>>> \bar c_0 \\
@VVV @VVV @VVV \\
p_1^*\mathcal{L}^{-1} \otimes e_1 \otimes p_1^*\pi_*R_1 @>>> 
A \otimes_{\mathbb{C}} (p_1^*\mathcal{L}^{-1} \otimes e_1) @>>> e_1
\end{CD}
\]
where 
\[
\begin{split}
&c_1 = (p_1^*\pi_*\pi^*\mathcal{L}^{-1} \otimes 
p_2^*p_{2*}(e_1 \otimes p_1^*\pi_*R_1))^{\Delta G} \\
&c_{10} = A_1 \otimes_{\mathbb{C}} 
(p_1^*\pi_*\pi^*\mathcal{L}^{-1} \otimes 
p_2^*p_{2*}(e_1 \otimes p_1^*\pi_*\pi^*\mathcal{O}_{\mathcal{X}}))^{\Delta G} 
\\
&c_0 = (p_1^*\pi_*\pi^*\mathcal{O}_{\mathcal{X}} \otimes 
p_2^*p_{2*}(e_1 \otimes p_1^*\pi_*\pi^*\mathcal{O}_{\mathcal{X}}))^{\Delta G} 
\\
&\bar c_1 = p_1^*\pi_*\pi^*\mathcal{L}^{-1} \otimes 
e_1 \otimes p_1^*\pi_*R_1 \\
&\bar c_{10} = A_1 \otimes_{\mathbb{C}} 
(p_1^*\pi_*\pi^*\mathcal{L}^{-1} \otimes 
e_1 \otimes p_1^*\pi_*\pi^*\mathcal{O}_{\mathcal{X}})
\\
&\bar c_0 = p_1^*\pi_*\pi^*\mathcal{O}_{\mathcal{X}} \otimes 
e_1 \otimes p_1^*\pi_*\pi^*\mathcal{O}_{\mathcal{X}}.
\end{split}
\]
We define a morphism 
$\bar \delta_0: c_0 \to e_1$ as the composition of the 
arrows in the last column.
Since the composition of the last row is zero, 
we have $\bar \delta_0 \delta_1 = 0$.

We have 
\[
\begin{split}
&\text{Hom}(\pi_*\pi^*\mathcal{L}^{-k} 
\boxtimes F(\pi_*R_k)[r], e_1) 
\cong \text{Hom}(p_2^*\Phi^{e_1}(\pi_*R_k)[r], 
e_1 \otimes p_1^*\pi_*\pi^*\mathcal{L}^k) \\
&\cong \text{Hom}(\Phi^{e_1}(\pi_*R_k)[r], 
\Phi^{e_1}(\pi_*\pi^*\mathcal{L}^k)) 
\cong \text{Hom}(\pi_*R_k[r], \pi_*\pi^*\mathcal{L}^k) \cong 0
\end{split}
\]
for any $r > 0$. 
If we apply the argument of the proof of Lemma~\ref{complex1lem} to a longer
complex
\[
0 \to c_m \to \dots \to c_0 \to e_1 \to 0
\]
we deduce that there exists a uniquely determined morphism 
$\epsilon: e' \to e_1$ such that $\epsilon \delta_0 = \bar \delta_0$.
Composing with the natural morphism $e \to e'$, 
we obtain a morphism $e \to e_1$.

Let $c$ be the cone of this morphism and assume that $c \not\cong 0$:
\[
e \to e_1 \to c \to e[1].
\]
We use the notation of Lemma~\ref{e-split}, where 
we assume additionally that $p_{2*}(c \otimes p_1^*\mathcal{A}) \ne 0$.
If we apply the functor $p_{2*}(\bullet \otimes p_1^*\mathcal{A})$ to the 
complex $c_{\bullet}$ with the morphism $\bar \delta_0: c_0 \to e_1$, 
then we obtain a complex $F(a_{\bullet})$ with a morphism
$\bar \delta_{0*}: F(a_0) \to \Phi^{e_1}(\mathcal{A})$, 
and a commutative diagram
\[
\begin{CD}
F(a_0) @>{\delta_{0*}}>> F(U(\mathcal{A})[m] \oplus \mathcal{A}) \cong 
\Phi^{e'}(\mathcal{A}) \\
@V=VV @V{\epsilon(\mathcal{A})}VV \\
F(a_0) @>{\bar \delta_{0*}}>> F(\mathcal{A}) \cong \Phi^{e_1}(\mathcal{A}).
\end{CD}
\]
Therefore, the induced morphism 
$p_{2*}(e \otimes p_1^*\mathcal{A}) \to p_{2*}(e_1 \otimes p_1^*\mathcal{A})$
is an isomorphism, hence $p_{2*}(c \otimes p_1^*\mathcal{A}) \cong 0$, 
a contradiction.

This is the end of the proof of Theorem~\ref{main}.


\section{Applications}

The following theorem is a generalization of \cite{DK}~Theorem~2.3 to the 
case of varieties with quotient singularities.

\begin{Thm}\label{application1}
Let $X$ and $Y$ be normal projective varieties with only quotient 
singularities, and $\mathcal{X}$ and $\mathcal{Y}$ the associated smooth 
stacks.
Assume that the bounded derived categories of coherent sheaves are 
equivalent as triangulated categories: $D^b(\text{Coh}(\mathcal{X})) 
\cong D^b(\text{Coh}(\mathcal{Y}))$.
Then the following hold:

(0) $\dim X = \dim Y$. 

(1) If $K_X$ (resp. $-K_X$) is nef, then $K_Y$ (resp. $-K_Y$) is also nef, 
and an equality on the numerical Kodaira dimension
$\nu(X) = \nu(Y)$ (resp. $\nu(X, -K_X) = \nu(Y, -K_Y)$) holds.

(2) If $\kappa(X) = \dim X$, i.e., $X$ is of general type, or if 
$\kappa(X, -K_X) = \dim X$, then $X$ and $Y$ are birationally equivalent. 
Moreover, there exist birational morphisms 
$f: Z \to X$ and $g: Z \to Y$ from a 
smooth projective variety $Z$ such that the canonical divisors are 
${\mathbb{Q}}$-linearly equivalent: $f^*K_X \sim_{\mathbb{Q}} g^*K_Y$.
\end{Thm}

\begin{proof}
The proof is parallel to that of \cite{DK}~Theorem~2.3, and
we only explain the outline.

By Theorem~\ref{main}, there exists an object $e \in 
D^b(\text{Coh}(\mathcal{X} \times \mathcal{Y}))$ 
such that $\Phi = \Phi^e_{\mathcal{X} \to \mathcal{Y}}: 
D^b(\text{Coh}(\mathcal{X})) 
\to D^b(\text{Coh}(\mathcal{Y}))$ is an equivalence.
Considering the right and left adjoint functors, 
we obtain an isomorphism of objects
\[
e^{\vee} \otimes p_1^*\omega_{\mathcal{X}}[\dim X]
\cong e^{\vee} \otimes p_2^*\omega_{\mathcal{Y}}[\dim Y].
\]
Part (0) follows immediately.

Let $\Gamma$ be the union of the 
supports of the cohomology sheaves $H^i(e^{\vee})$ for all $i$, 
and let $\hat{\nu}_1 : \tilde{\mathcal{Z}}_1 \to \mathcal{Z}_1$ 
be the normalization of an 
irreducible component of $\Gamma$ which is at the same time 
an irreducible component of the support of some $H^i(e^{\vee})$.
Since $\Phi^e_{\mathcal{X} \to \mathcal{Y}}$ is an equivalence, 
there exists $\mathcal{Z}_1$  such that
the projection $\mathcal{Z}_1 \to \mathcal{Y}$ is surjective and that 
\[
\hat{\nu}_1^*p_1^*\omega_{\mathcal{X}}^{\otimes m_1} \cong 
\hat{\nu}_1^*p_2^*\omega_{\mathcal{Y}}^{\otimes m_1}
\]
where $m_1$ is the rank of $\hat{\nu}_1^*H^i(e^{\vee})$.
Let $Z_1$ be the image of $\mathcal{Z}_1$ on $X \times Y$ and
$\nu_1 : \tilde Z_1 \to Z_1$ the normalization.

If $K_X$ is nef or anti-nef, then 
$\nu_1^*p_1^*K_X \sim_{\mathbb{Q}} \nu_1^*p_2^*K_Y$ 
is also nef or anti-nef, hence so is $K_Y$.  
We have also 
$\nu(X, \pm K_X) \ge \nu(\tilde Z_1, \pm \nu_1^*p_2^*K_Y) = \nu(Y, \pm K_Y)$, 
thus $\nu(X, \pm K_X) = \nu(Y, \pm K_Y)$ by symmetry.

If $\kappa(X, \pm K_X) = \dim X$, then there exist an ample 
$\mathbb{Q}$-Cartier divisor $A$ 
and an effective $\mathbb{Q}$-Cartier divisor $B$ on $X$ 
such that $\pm K_X \sim_{\mathbb{Q}} A + B$ by Kodaira's lemma.
We take $\mathcal{Z}_1$ which dominates $\mathcal{X}$ instead of $\mathcal{Y}$.
Then the projection 
$p_2 \vert_{Z_1}: Z_1 \to Y$ 
is quasi-finite on $Z_1 \setminus p_1^{-1}(\text{Supp}(B))$.
It follows that $\dim Z_1 = \dim X$ and 
$Z_1$ also dominates $Y$.
The set $\Gamma \cap p_1^{-1}(x)$ consisits
of a single point for a general point $x \in X$, 
and $Z_1$ becomes a graph of a birational map.
If we take $Z$ to be its resolution of singularities, 
then the conclusion holds.
\end{proof}

We consider a generalization of a result in \cite{BO} on the autoequivalence 
group:

\begin{Thm}\label{application2}
Let $X$ be a normal projective variety with only quotient 
singularities, and $\mathcal{X}$ the associated smooth stack.
Assume that $K_X$ or $-K_X$ is ample and that $K_X$ generates the local 
divisor class group at each point $x \in X$.
Then the group of autoequivalences 
$\text{Autoeq}(D^b(\text{Coh}(\mathcal{X})))$ is 
isomorphic to the semidirect product of the automorphism group $\text{Aut}(X)$
and the trivial factor $\text{Pic}(\mathcal{X}) \times \mathbb{Z}$.
\end{Thm}

\begin{proof}
The theorem can be proved in the same way as in \cite{BO}.
But we present here an alternative proof using Theorem~\ref{main}
as in \cite{DK}~Remark~2.4~(1).

We use the notation of the proof of Theorem~\ref{application1}.
Since we can take $B=0$, $Z_1$ becomes a graph of an automorphism of $X$, 
which induces an automorphism $h$ of $\mathcal{X}$.
We have that 
\begin{equation}\label{Serre}
\Phi(a(jK_{\mathcal{X}})) \cong 
\Phi(a)(jK_{\mathcal{X}})
\end{equation}
for any object $a$ and integer $j$.

We claim that the support of $e$ coincides with the graph of $h$.
Indeed, if not, then there exists a point $x \in \mathcal{X}$ such that 
the support of $e \otimes p_1^*\mathcal{O}_x$ is not connected.
Since $K_X$ generates the local class groups, it follows that 
the support of $\Phi(\mathcal{O}_x)$ is not connected by 
$(\ref{Serre})$, a contradiction.

Let $x \in \mathcal{X}$ be any point, $r_x$ the index of $K_X$ at $\sigma(x)$,
and $j$ an integer such that $0 \le j < r_x$.
Then $\Phi(\mathcal{O}_x)$ is a shift of a sheaf, because we have otherwise
$\text{Hom}_{D^b(\mathcal{X})}^p(\Phi(\mathcal{O}_x), 
\Phi(\mathcal{O}_x(jK_X))) \ne 0$ for some $p < 0$ and $j$.
Moreover, since 
\[
\text{Hom}_{D^b(\mathcal{X})}(\Phi(\mathcal{O}_x), \Phi(\mathcal{O}_x(jK_X))) 
\cong \begin{cases} \mathbb{C} &\text{ if } j = 0 \\
0 &\text{ if } 0 < j < r_x
\end{cases}
\]
$\Phi(\mathcal{O}_x)$ is a shift of a skyscraper sheaf of length $1$.

Assume that $e$ is not a shift of a sheaf, and let 
$i_1$ be the maximum of the integers $i$
such that $H^i(e) \ne 0$.
We consider a spectral sequence
\[
E_2^{p,q} = \text{Tor}_{-p}(H^q(e), p_1^*\mathcal{O}_x) 
\Rightarrow e \otimes p_1^*\mathcal{O}_x.
\]
Take a general point $x$ of the support of $p_{1*}H^{i_1}(e)$ 
which is also contained in the support of $p_{1*}H^{i_2}(e)$ for 
some $i_2 < i_1$.
We assume that $i_2$ is the largest such integer.
If the support of $p_{1*}H^{i_1}(e)$ is a proper subvariety of $\mathcal{X}$, 
then both $E_2^{0,i_1}$ and $E_2^{-1,i_1}$ do not vanish, and they survive at 
$E_{\infty}$.
Otherwise, both $E_2^{0,i_1}$ and $E_2^{0,i_2}$ 
do not vanish and survive at $E_{\infty}$.
In any case, we conclude that $\Phi(\mathcal{O}_x)$ is not 
a shift of a skyscraper sheaf of length $1$, a contradiction.
Hence $e[-i_1]$ is a sheaf.
Moreover, there exists an invertible sheaf $\mathcal{M}$ 
on $\mathcal{X}$ such that
$\Phi^e(a) \cong h_*(a) \otimes \mathcal{M}[i_1]$ for any 
$a \in D^b(\mathcal{X})$.
\end{proof}

\begin{Rem}
Let $Y$ be a minimal algebraic surface of general type, and 
$X$ its canonical model.

(1) By \cite{BKR}, the derived categories $D^b(\text{Coh}(\mathcal{X}))$ 
and $D^b(\text{Coh}(Y))$ 
are equivalent.
Though $K_X$ is ample, $X$ and $Y$ may not be isomorphic.
This is a difference from the case of smooth varieties (\cite{BO}); 
if $X$ is smooth and $\pm K_X$ is ample, then any variety 
$Y$ with an equivalent derived category is isomorphic to $X$.
The reason is that, 
there is no crepant birational map, which is not an isomorphism,  
from a smooth projective variety
whose canonical or anti-canonical divisor is ample.
 
(2) The condition on the local class groups in Theorem~\ref{application2}
cannot be removed. 
For example, assume that $X$ has an ordinary double point $x_0 \in X$.
Let $\mathcal{O}_{x_0}(-1)$ be a skyscraper stacky sheaf over $x_0$ 
which has a non-trivial action of the stabilizer group, and
$\mathcal{O}_{\Delta \mathcal{X}}$ the object on the product 
$\mathcal{X} \times \mathcal{X}$ corresponding to the 
identity functor.
Then the object
\[
e = \text{cone}\{\mathcal{O}_{x_0}(-1)^{\vee} \boxtimes \mathcal{O}_{x_0}(-1) 
\to \mathcal{O}_{\Delta \mathcal{X}}\}
\]
corresponds to the twisting $T \in \text{Autoeq}(D^b(\text{Coh}(\mathcal{X})))$ 
defined by
\[
T: a \mapsto \text{cone}\{\text{Hom}(\mathcal{O}_{x_0}(-1), a) 
\otimes \mathcal{O}_{x_0}(-1) \to a\}
\]
for $a \in D^b(\text{Coh}(\mathcal{X}))$ (\cite{ST}).
We note that 
\[
\mathcal{O}_{x_0}(-1)^{\vee} \cong \mathcal{O}_{x_0}(-1)[-2].
\]
\end{Rem}

\begin{Rem}
The main results of this paper hold without change for more general 
situation where the stabilizer groups of the smooth stacks have fixed loci of 
codimension $1$.
This extended version will be used in a forthcoming paper.
\end{Rem}


Department of Mathematical Sciences, University of Tokyo, 

Komaba, Meguro, Tokyo, 153-8914, Japan 

kawamata@ms.u-tokyo.ac.jp

\end{document}